\documentclass[11pt]{amsart}
\usepackage[normalem]{ulem}%\uline,\uuline,\uwave,\sout,\xout
\usepackage[usenames]{color}
\usepackage{amssymb}
\usepackage{graphicx}
\usepackage{url}
\usepackage{amsrefs}
\usepackage{tikz}
\usetikzlibrary{arrows.meta}
\usetikzlibrary{decorations.markings}

 \newcommand{\R}{\mathbb R}
 
\newcommand{\ep}{\varepsilon}

\theoremstyle{plain} \newtheorem{thm}{Theorem}
 \newtheorem{prop}[thm]{Proposition}

\theoremstyle{definition} \newtheorem{defn}[thm]{Definition}

\newtheorem{ex}[thm]{Example} 

\theoremstyle{remark} 

\title{A generalization of the Brouwer plane translation theorem}

\author{Jim Wiseman}
\address{Department of Mathematics \\ Agnes Scott College \\ Decatur, GA 30030}
\email{jwiseman@agnesscott.edu}

\begin{document}

\begin{abstract}
We show that if an orientation-preserving homeomorphism of the plane has a topologically chain recurrent point, then it has a fixed point, generalizing the Brouwer plane translation theorem.
\end{abstract}

\maketitle

Let $f:\R^2\to \R^2$ be an orientation-preserving homeomorphism.  The Brouwer plane translation theorem  (\cite{Brouwer}) says that if $f$ has a periodic point, then it has a fixed point.  Fathi weakened the condition to having a non-wandering point (\cite{Fathi}).  Mai, Yang, and Zen weakened it to having a BP-chain recurrent point
(Definition~\ref{defn:BPCR}), and showed that just having a chain recurrent point was not sufficient (\cite{MYZ}).  (See \cites{Franks92,Brown84,Slaminka,Guillou,MYZ,LeCalvez} for other proofs and more on the history of the theorem.)
We weaken the condition further by proving that if $f$ has a topologically chain recurrent point, then it has a fixed point. 

We begin with definitions.  For the remainder of the paper, let $d_0$ be the standard metric on $\R^2$.

Let $(X,d)$ be a metric space and $f:X\to X$ a continuous map.  For $\ep>0$, an \emph{$\ep$-chain} from the point $p$ to the point $q$ is a sequence $(p_0=p, p_1, \ldots, p_n=q)$ such that $d(f(p_{i-1}), p_i) < \ep$ for $1\le i \le n$.  A point $p$ is \emph{chain recurrent} if for every $\ep>0$ there is an $\ep$-chain from $p$ to itself.

Mai, Yang, and Zen give an example of an orientation-preserving homeomorphism of the plane where every point is chain recurrent  but there is no fixed point  (\cite{MYZ}*{Example~3.1}).  We give a simpler example here.

\begin{ex} \label{ex:translation}
Consider the orientation-preserving homeomorphism $f(x,y) = (x+e^y, y)$.  There are clearly no fixed points.  However, every point is chain recurrent:  for $\ep>0$, we can create an $\ep$-chain from any point back to itself by moving right and down ($p_i = f(p_{i-1}) - (0,\ep)$) until $e^y < \ep$, then sufficiently far to the left ($p_i = f(p_{i-1}) - (\ep,0)$), and then right and up until we return to the original point.
\end{ex}

If $X$ is compact, then chain recurrence is independent of the choice of metric (compatible with the topology) (see, for example, \cite{Franks87}).  This is not true on noncompact spaces, as the following example shows.

\begin{ex} \label{ex:circ}
Let $f(x) = x+1$ be translation on the real line $\R$ with the usual metric.  Then no point is chain recurrent.  However, if we think of $\R$ as the unit circle with a point removed and give it the metric induced by the circle metric, then every point is chain recurrent.
\end{ex}

To avoid this dependence on the choice of metric, we use the notion of topological chain recurrence from \cite{DLRW}, which is equivalent to Hurley's strong chain recurrence from \cite{Hurley:1991}, as shown in \cite{DLRW}.

Let $\Delta = \{(p,p): p \in \R^2\}$ be the diagonal in the product $\R^2\times\R^2$.  For any neighborhood $N$ of $\Delta$, we define an \emph{$N$-chain} to be a sequence $(p_0, p_1, \ldots, p_n)$ of points in $\R^2$  ($n\ge1$) such that $(f(p_{i-1}), p_i) \in N$ for $i=1,\ldots,n$.  A point $p\in\R^2$ is \emph{topologically chain recurrent} if for every neighborhood $N$ of $\Delta$, there is an $N$-chain from $p$ to itself.

Note that if we take $N$ to be the $\ep$-neighborhood of $\Delta$, then an $N$-chain is equivalent to an $\ep$-chain.  Thus topological chain recurrence implies (ordinary) chain recurrence, but, as is shown in \cite{DLRW}, the converse is not true.

We can now state the main result.

\begin{thm} \label{thm:main}
Let $f:\R^2\to R^2$ be an orientation-preserving homeomorphism. If $f$ has a topologically chain recurrent point, then it has a fixed point.
\end{thm}

The proof relies on the following result of Franks.

\begin{defn}
A \emph{periodic disk chain} for $f$ is a finite set $U_0,U_2,\ldots,U_n = U_0$ of topological disks in $\R^2$ such that
\begin{enumerate}
\item $f(U_i) \cap U_i = \emptyset$ for $i=0,\ldots,n$.
\item  If $i\ne j$ then either $U_i=U_j$ or $U_i\cap U_j = \emptyset$.
\item For $i=1,\ldots,n$, there exists $m_{i-1} >0$ with $f^{m_{i-1}}(U_{i-1})\cap U_{i} \ne \emptyset$.
\end{enumerate}
\end{defn}

\begin{prop}[\cite{Franks}*{Proposition~1.3}] \label{prop:Franks}  
Let $f:\R^2\to \R^2$ be an orientation-preserving homeomorphism. If $f$ possesses a periodic disk chain, then it has a fixed point.
\end{prop}

\begin{proof}[Proof of Theorem~\ref{thm:main}]
For a neighborhood $N$ of the diagonal $\Delta$ and a point $p\in\R^2$, define the set $N(p)$ in $\R^2$ by $N(p) = \{q: (p,q) \in N\}$.

Assume that $f$ has no fixed point.  Then there exists a neighborhood $N$ of $\Delta$ such that $N(p) \cap f(N(p)) = \emptyset$ for all $p \in \R^2$.

  Let $p'$ be a topologically chain recurrent point.  Then there is an $N$-chain $(p'=p_0, p_1, \ldots, p_n=p')$ from $p'$ to itself.  Define $U_i = N(f(p_{i-1}))$, $1\le i \le n$, and $U_0=U_n = N(f(p_{n-1}))$.
  After some modification, if necessary, the set $U_0,U_2,\ldots,U_n = U_0$ will give a periodic disk chain, and the result will follow from Proposition~\ref{prop:Franks}.  The details are as follows.

Property (1) of the definition of periodic disk chain will be true if we use $N$ or any subset of $N$.  
We need to ensure that each $U_i$ is a topological disk.  If that is not the case, we can proceed as follows (as in \cite{DLRW}*{Proposition 34(3)}).  We will construct a continuous function $\delta:\R^2 \to (0,\infty)$ such that $B_\delta \subset N$, where $B_\delta = \{(p,q) : d(p,q) < \delta(p)$; then, if necessary, we can replace $N$ by $B_\delta$, and each $U_i$ will be a topological disk.  First, we can assume that $N(p) \subsetneq \R^2$ for all $p$, and  define a function $h:\R^2\to(0,\infty)$ by $h(x) = d(x, \R^2-N(x))$.  Then, since $h$ is lower semicontinuous, the function $\delta(p)=\inf\{h(q) + d(p,q): q\in \R^2\}$ is continuous (\cite{Ziemer}, Theorem~3.66), with $0<\delta(p) < h(p)$ for all $p$; thus $B_\delta \subset B_h \subset N$.

To get property (2), we may need to further modify $N$ or the $N$-chain.  
Assume that $U_i \cap U_j \ne \emptyset$ for some $i<j$ (other than $U_0=U_n$).  We consider several cases.
If neither $p_j$ nor $f(p_{j-1})$ is in $U_i \cap U_j$, then we can shrink $N$ slightly so that $U_i$ and  $U_j$ no longer intersect but we still have an $N$-chain.  If $p_j \in U_i \cap U_j$, then we can remove a segment of the middle of the $N$-chain, to eliminate the intersection, and $(p'=p_0,p_1,\ldots, p_{i-1},p_j,\ldots, p_n=p')$ will still be an $N$-chain.  
After resolving all of the intersections in the previous cases, there may be no more intersections, in which case we have a periodic disk chain, with all of the $m$'s equal to 1 for property (3) of the definition.  Otherwise, we  consider finally the case $f(p_{j-1}) \in U_i \cap U_j$.  Take the maximum such $i$,  and then take the minimum such $j>i$.  Then $U_i, U_{i+1},\ldots, U_{j-1}, U_i$ is a periodic disk chain, again with all of the $m$'s equal to 1 for property (3).

%Assume that $N(p_i) \cap N(p_j) \ne \emptyset$ for some $i<j$ (other than $N(p_0)=N(p_n)$.  We consider several cases.  If neither $p_i$ nor $f(p_{i-1})$ is in $N(p_j)$, then we can shrink $N$ slightly so that $N(p_i)$ and  $N(p_j) $ no longer intersect but we still have an $N$-chain.  If $f(p_{i-1})\in N(p_j)$, then we can remove a segment of the middle of the $N$-chain, to eliminate the intersection, and $(p'=p_0,p_1,\ldots, p_{i-1},p_j,\ldots, p_n=p')$ will still be an $N$-chain.  The last case is $p_i\in N(p_j)$.  Again, we consider cases.  If neither $p_j$ nor $f(p_{j-1})$ is in $N(p_i)$, then again we can shrink $N$ slightly so that $N(p_i)$ and  $N(p_j) $ no longer intersect but we still have an $N$-chain.  If $p_j \in N(p_i)$, then again we can remove a segment of the middle of the $N$-chain, to eliminate the intersection, and $(p'=p_0,p_1,\ldots, p_{i-1},p_j,\ldots, p_n=p')$ will still be an $N$-chain.  After resolving all of the intersections in the previous cases, there may be no more intersections, in which case we have a periodic disk chain, with all of the $m$'s equal to 1 for property (3) of the definition.  Otherwise, we  consider finally the case when $f(p_{j-1}) \in N(p_i)$.  For the given $i$,  take the minimum such $j$.  Then $N(p_i), N(p_{i+1}),\ldots, N(p_{j-1}), N(p_i)$ is a periodic disk chain, again with all of the $m$'s equal to 1 for property (3).

Since we have a periodic disk chain, Proposition~\ref{prop:Franks} implies that $f$ has a fixed point.
  
\end{proof}

In \cite{MYZ}, Mai, Yang, and Zen prove the following result.

\begin{defn} \label{defn:BPCR}
Let $(X,d)$ be a metric space and $f:X\to X$ a continuous map. 
A point $p$ is \emph{bounded-perturbation (BP) chain recurrent} if there exists a bounded set $W$ such that for every $\ep>0$, there is an $\ep$-chain from $p$ to $p$ with the property that the only jumps in the chain are at points in $W$, that is, $p_{i} = f(p_{i-1})$ if $f(p_{i-1})\not \in W$.
\end{defn}

\begin{thm}[\cite{MYZ}*{Theorem~3.5}] \label{thm:BPCRresult}
Let $f:\R^2\to R^2$ be an orientation-preserving homeomorphism. If $f$ has a  BP-chain recurrent point, then it has a fixed point.
\end{thm}

We will show that Theorem~\ref{thm:main} is stronger than Theorem ~\ref{thm:BPCRresult}, in the sense that bounded-perturbation chain recurrence implies topological chain recurrence, but not vice versa.

The definition of BP-chain recurrence, and the result in Theorem ~\ref{thm:BPCRresult}, depend on the metric.  
See Example~\ref{ex:circ} -- with the usual metric, no point is BP-chain recurrent, but with the metric induced from the circle metric, every point is BP-chain recurrent.
As an example on $\R^2$, if we take the bounded metric $\overline d$ given by $\overline d(p,q) = \max(d_0(p,q), 1)$, then every set is bounded and so BP-chain recurrence is equivalent to chain recurrence.  Example~\ref{ex:translation} gives an orientation-preserving homeomorphism of the plane with  no fixed point.  With the usual metric $d_0$, it has no BP-chain recurrent points, but with the bounded metric $\overline d$ every point is BP-chain recurrent.

To eliminate this dependence on the metric, we make the following definition, which simply replaces ``bounded'' with ``compact'' and  is thus equivalent to BP-chain recurrence for the standard metric on $\R^2$.

\begin{defn}
Let $(X,d)$ be a metric space and $f:X\to X$ a continuous map. 
A point $p$ is \emph{compact-perturbation chain recurrent} if there exists a compact set $W$ such that for every $\ep>0$, there is an $\ep$-chain from $p$ to $p$ with the property that the only jumps in the chain are at points in $W$, that is, $p_{i} = f(p_{i-1})$ if $f(p_{i-1})\not \in W$.  \end{defn}

\begin{prop}
Let $X$ be a metric space.  Then compact-perturbation chain recurrence is independent of the choice of metric (compatible with the topology).
\end{prop}

\begin{proof}
Let $d_1$ and $d_2$ be two metrics compatible with the topology.  The identity map sending $(X,d_1)$ to $(X,d_2)$ is uniformly continuous on any compact set $W$. Thus for any $\ep>0$ there exists a $\delta$ such that if $p,q \in W$ and $d_1(p,q) < \delta$, then $d_2(p,q) < \epsilon$.  Therefore a point is compact-perturbation  chain recurrent for $d_2$ if it is compact-perturbation  chain recurrent for $d_1$; the converse follows by symmetry.

\end{proof}

\begin{thm} \label{thm:CPCRimpliesTCR}
Let $X$ be a metric space.
Compact-perturbation chain recurrence implies topological chain recurrence, but not vice versa.
\end{thm}

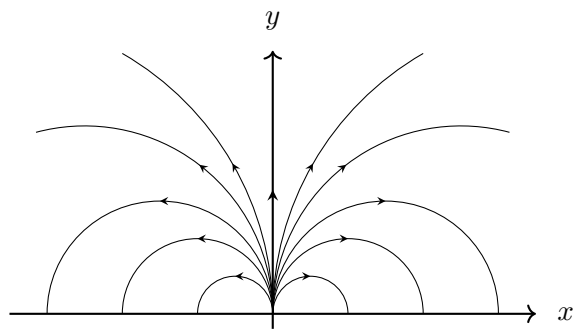
\begin{figure} 
\begin{tikzpicture}
[decoration={markings, 
	mark= at position 0.5 with {\arrow{stealth}}}
] 
   
    % Draw axes
        \draw [->,thick,postaction={decorate}] (0,-.2)   -- (0,3.5) node [label=above:{$y$}] {}; 
  	\draw [->,thick] (-3.5,0)   -- (3.5,0) node [label=right:{$x$}] {}; 
  	 \draw [postaction={decorate}] (0,0) arc (180:0:.5);
	\draw [postaction={decorate}] (0,0) arc (180:0:1);
	\draw [postaction={decorate}] (0,0) arc (180:0:1.5);
	\draw [postaction={decorate}] (0,0) arc (180:75:2.5);
	\draw [postaction={decorate}] (0,0) arc (180:120:4);
	\draw [postaction={decorate}] (0,0) arc (0:180:.5);
	\draw [postaction={decorate}] (0,0) arc (0:180:1);
	\draw [postaction={decorate}] (0,0) arc (0:180:1.5);
	\draw [postaction={decorate}] (0,0) arc (0:105:2.5);
	\draw [postaction={decorate}] (0,0) arc (0:60:4);
\end{tikzpicture}

\caption{The flow in the upper half-plane for the time-one map in Example~\ref{ex:CRnotCPCR}; points in the lower half-plane are fixed.  Points on the positive $y$-axis are topologically chain recurrent but not compact-perturbation chain recurrent.} \label{fig:TCRnotCPCR}
\end{figure}

\begin{proof}
Assume that $p$ is compact-perturbation chain recurrent for a map $f$.  Let $W$ be a compact set such that for every $\ep>0$, there is an $\ep$-chain from $p$ to $p$ with the property that the only jumps in the chain are at points in $W$.  Let $N$ be any neighborhood of $\Delta$.  Since $W$ is compact, there exists a $\delta>0$ such that the $\delta$-ball around $q$ is contained in $N(q)$ for every $q\in W$.  There is a $\delta$-chain from $p$ to $p$ with the only jumps at points in $W$; this chain is also an $N$-chain, and since $N$ was arbitrary, $p$ is topologically chain recurrent.

The following example shows that topological chain recurrence does not imply compact-perturbation chain recurrence.

\begin{ex} \label{ex:CRnotCPCR}
Let $f:\R^2\to\R^2$ be the orientation-preserving homeomorphism of the plane defined as the time-one map of the flow shown in Figure~\ref{fig:TCRnotCPCR}, which is the identity in the lower half-plane $\{y\le0\}$.  Any point on the positive $y$-axis is topologically chain recurrent: make a small jump to one of the semi-circles, follow the semi-circle to the $x$-axis, make small jumps along the $x$-axis back to the origin, then jump to the positive $y$-axis and follow it back to the starting point.  However, no point on the positive $y$-axis is compact-perturbation chain recurrent: as $\ep \to 0$, the jumps are limited to semi-circles of larger and larger radius, which then require jumps further and further out the $x$-axis in order to return to the origin.
\end{ex}

\end{proof}

Thus, by Theorem~\ref{thm:CPCRimpliesTCR}, Theorem ~\ref{thm:BPCRresult} is a special case of Theorem~\ref{thm:main}.

% \bib, bibdiv, biblist are defined by the amsrefs package.

\end{document}